# AUTOMORPHISMS OF THE ENDOMORPHISM SEMIGROUP OF A FREE COMMUTATIVE ALGEBRA

A. BELOV-KANEL[1], R. LIPYANSKI[2]

ABSTRACT. We describe the automorphism group of the endomorphism semigroup $\mathrm{End}(K[x_1, \ldots, x_n])$ of ring $K[x_1, \ldots, x_n]$ of polynomials over an *arbitrary* field $K$. A similar result is obtained for automorphism group of the category of finitely generated free commutative-associative algebras of the variety $\mathcal{CA}$ commutative algebras. This solves two problems posed by B. Plotkin ( [18], Problems 12 and 15).

More precisely, we prove that if $\varphi \in \mathrm{AutEnd}(K[x_1, \ldots, x_n])$ then there exists a semi-linear automorphism $s : K[x_1, \ldots, x_n] \to K[x_1, \ldots, x_n]$ such that $\varphi(g) = s \circ g \circ s^{-1}$ for any $g \in \mathrm{End}(K[x_1, \ldots, x_n])$. This extends the result by A. Berzins obtained for an infinite field $K$.

## 1. INTRODUCTION

We describe the group $G = \mathrm{Aut}(\mathrm{End}(K[x_1, \ldots, x_n])$, where $K$ is an arbitrary field. A similar result is obtained also for automorphism group of the category of finitely generated free commutative-associative algebras of the variety commutative algebras. This solves two problems posed by B. Plotkin ( [18], Problems 12 and 15).

More precisely, we prove that if $\varphi \in \mathrm{AutEnd}(K[x_1, \ldots, x_n])$ then there exists a semi-linear automorphism $s : K[x_1, \ldots, x_n] \to K[x_1, \ldots, x_n]$ such that $\varphi(g) = s \circ g \circ s^{-1}$ for any $g \in \mathrm{End}(K[x_1, \ldots, x_n])$ (see Theorem 3.8). Here "semi-linearity" means that $s$ is a composition of an automorphism of the field $K$ and an automorphism of the ring $K[x_1, \ldots, x_n]$. We note that for an infinite ground field $K$ is infinite such result was obtained earlier by A. Berzins [3].

A problem of description of the group $G = \mathrm{Aut}(\mathrm{End}(K[x_1, \ldots, x_n])$ is also interesting in the context of Universal Algebraic Geometry (UAG). Let $\Theta$ be a variety of algebras over a field $K$ and $F = F(X)$ be a free algebra from $\Theta$ generated by a finite subset $X$ of some infinite universum $X^0$. We refer to [17, 18] (see also [8]) for the Universal Algebraic Geometry (UAG) notions used in our work.

If an algebra $G$ belongs to $\Theta$ one can consider the category of algebraic sets $K_\Theta(G)$ over $G$. Objects of this category are algebraic sets in affine space over $G$; the category $K_\Theta(G)$ defines a geometry of the algebra $G$ in $\Theta$. One of the main problems in UAG is to determine whether two different algebras $G_1$ and $G_2$ have the same geometry. The coincidence of geometries means that the categories $K_\Theta(G_1)$ and $K_\Theta(G_2)$ are equivalent. It is known that coincidence of geometries of $G_1$ and $G_2$ is determined by the structure of the group $\mathrm{Aut}\,\Theta^0$, where $\Theta^0$ is the category of free finitely generated algebras of $\Theta$. On the other hand, there is a natural relation between the structure of the groups $\operatorname{Aut}\operatorname{End} F$ and $\operatorname{Aut}\Theta^0$. The structure of the latter is determined by the group $\operatorname{Aut}\operatorname{End} F$. Should be mentioned that a problem of investigation of the groups $\operatorname{Aut}\operatorname{End} F$, $F \in \Theta$, for different varieties $\Theta$ is quite interesting by itself and has been considered in many papers (see [1]-[3], [5], [8]-[11], [13]-[19],[22]

Let $\mathcal{CA}$ be the variety of a commutative-associative algebras with 1 over a field $K$, $A = K[x_1, ..., x_n]$ be a free commutative-associative algebra in $\mathcal{CA}$ freely generated over $K$ by a set $X = \{x_1, ..., x_n\}$, i.e., a polynomial algebra in variables $x_1, ..., x_n$. In this work we obtain a description of the group $Aut\,\mathcal{CA}^0$ of automorphisms of the category $\mathcal{CA}$. A similar result for a polynomial algebra $A$ over an infinite field $K$ was also obtained earlier in [3].

Our description is based on new characteristics of endomorphisms of $A$ such as *rank* of endomorphisms of $A$. We discuss external and internal definition of this notation. The former are expressed in terms of the action of the semigroup $\operatorname{End} A$ on $A$, while the latter can be written in terms of the semigroup itself. This approach allows us to describe the above mentioned properties of endomorphisms of $A$ in an invariant manner and paves the way for proof of the main assertions in the paper: the group $\operatorname{Aut}\operatorname{End} A$ is generated by semi-inner of $\operatorname{End} A$.

Our approach employs this technique (developed in [5, 9]) supplemented by algebro-geometric methods of investigations

## 2. On the endomorphism semigroup of a free associative-commutative algebra

**2.1. Rank of an endomorphism of polynomial algebra.** Let $A = K[x_1, \ldots, x_n]$ be a free commutative-associative algebra over a field $K$ generated by $X = \{x_1, \ldots, x_n\}$ (below *polynomial algebra* over $K$ in variables $X$). Earlier, in [5], we defined *the endomorphism* of free associative algebra $K\langle x_1, \ldots, x_n\rangle$ *of rank* 0 *and* 1. In this section we introduce a definition of *endomorphisms of arbitrary rank $m$* in a free commutative-associative $K[x_1, \ldots, x_n]$.

First, we introduce the "external" and "internal" definitions of *rank* of endomorphism $\varphi$ of algebra $A$ and show their equivalence.

**Definition 2.1.** ("External" definition of an endomorphism of rank $m$.) An endomorphism
$$\varphi : A \to A$$
has **rank** $m$ if $\operatorname{trdeg}(\operatorname{Im}\varphi) = m$, i.e., the transcendence degree of the $K$-algebra $M = \operatorname{Im}\varphi \subseteq A$ is equal to $m$. We denote this as $\operatorname{rk}(\varphi) = m$. It is evident that there exist endomorphisms of $K[x_1, \ldots, x_n]$ of arbitrary rank $\leq n$. For instance, the identical mapping on $K[x_1, \ldots, x_n]$ is the endomorphism of rank $n$.

For the internal definition of rank $m$ endomorphisms, we need to define a congruence on the semigroup $\operatorname{End}(A)$ with respect to a fixed endomorphism $\varphi$ of $A$.

**Definition 2.2.** Endomorphisms $\varphi_1$ and $\varphi_2$ of $A$ are $\varphi$-*equivalent* if $\varphi\varphi_1 = \varphi\varphi_2$. In this case we write $\varphi_1 \backsim_\varphi \varphi_2$.

It is clear that $\backsim_\varphi$ is an equivalence relation on $\operatorname{End} A$. Let $S$ be the set of all $\varphi$-equivalences on $\operatorname{End} A$. We determine the preorder $\trianglelefteq$ on the set $S$ as follows. We say that $\backsim_\phi \trianglelefteq \backsim_\psi$, where $\phi, \psi \in \operatorname{End} A$, if

$$\phi\varphi_1 = \phi\varphi_2 \Rightarrow \psi\varphi_1 = \psi\varphi_2,$$

for any $\varphi_1, \varphi_2 \in \operatorname{End} A$. The preorder $\trianglelefteq$ can be extended up to the order $\preceq$ on the quotient set $\tilde{S} = S/R$ under equivalence $R$, where $\backsim_\phi R \backsim_\psi$ if and only if $\backsim_\phi \trianglelefteq \backsim_\psi$ and $\backsim_\psi \trianglelefteq \backsim_\phi$. Denote by $\backsim_{\psi_R}$ the $R$-equivalence class of a relation $\backsim_\psi$.

**Definition 2.3.** We say that $\phi \preceq \psi$ iff $\backsim_{\phi_R} \preceq \backsim_{\psi_R}$.

**Definition 2.4.** We say that $\phi \prec \psi$ if $\backsim_{\phi_R} \preceq \backsim_{\psi_R}$ and $\backsim_{\psi_R} \not\preceq \backsim_{\phi_R}$.

It is clear that relations $\preceq$ and $\prec$ are an order and a strong order, respectively, on $\operatorname{End} A$. Note that the smaller endomorphism $\varphi$ (in the sense of $\preceq$) corresponds to stronger equivalence relation $\sim_\varphi$. The proof of the following Lemma is straightforward.

**Lemma 2.5.** *Let $\varphi = (\varphi_1(\vec{x}) \ldots, \varphi_n(\vec{x}))$ and $\phi = (\psi_1(\vec{x}), \ldots, \psi_n(\vec{x}))$ be two endomorphisms of $K[x_1, \ldots, x_n]$. Then*
  (1) *$\phi \sim \psi$ iff for all $H(\vec{x}) \in K[x_1, \ldots, x_n]$ the condition $H(\varphi_1(\vec{x}) \ldots, \varphi_n(\vec{x})) = 0$ is equivalent to $H(\psi_1(\vec{x}), \ldots, \psi_n(\vec{x})) = 0$.*
  (2) *$\phi \preceq \psi$ iff for all $H(\vec{x}) \in K[x_1, \ldots, x_n]$ the condition $H(\varphi_1(\vec{x}), \ldots, \varphi_n(\vec{x})) = 0$ implies $H(\psi_1(\vec{x}), \ldots, \psi_n(\vec{x})) = 0$.*
  (3) *$\phi \prec \psi$ iff for all $H(\vec{x}) \in K[x_1, \ldots, x_n]$ the condition $H(\varphi_1(\vec{x}) \ldots, \varphi_n(\vec{x})) = 0$ implies $H(\psi_1(\vec{x}), \ldots, \psi_n(\vec{x})) = 0$ and there exists $R(\vec{x}) \in K[x_1, \ldots, x_n]$ such that $R(\varphi_1(\vec{x}), \ldots, \varphi_n(\vec{x})) = 0$ but $H(\psi_1(\vec{x}), \ldots, \psi_n(\vec{x})) \neq 0$.*

**Definition 2.6.** ("Internal" definition of an endomorphism of rank $m$.) An endomorphism $\psi: A \to A$ is of *rank* $m$, if maximum of the lengths of all chains of endomorphisms of $A$ of the form

$$(2.1) \qquad \psi \precnsim \psi_{m-1} \precnsim \cdots \precnsim \psi_1 \precnsim \psi_0,$$

is equal to $m$. If there is no endomorphism $\psi$ such that $\psi \precnsim \psi_0$, then $\psi$ has *rank* 0.

**Remark 2.7.** If $\operatorname{rk}(\varphi) = 0$, then image of $\varphi$ is the ground field. The definition of endomorphisms of rank 0 and 1 for associative commutative algebra are in accordance with the definition for a free associative algebra given in [5]. The internal definition of rank 0 is pretty similar.

**Proposition 2.8.** *Definitions 2.6 and 2.1 are equivalent.*

We precede the proof of this proposition by several lemmas. Denote by $\mathbf{A}^n_K$ an n-dimensional affine space over the algebraic closure $\bar{K}$ of the field $K$. It is clear that $\mathbf{A}^n_K \simeq \operatorname{Specm}(K[x_1, \ldots, x_n])$, where $\operatorname{Specm}(K[x_1, \ldots, x_n]$ is the set of all maximal ideals. Let us investigate the algebro-geometric properties of polynomial endomorphisms of $K[x_1, \ldots, x_n]$ and their relation to polynomial maps of $\mathbf{A}^n_K$ into itself.

Each endomorphism $\varphi: K[x_1, \ldots, x_n] \to K[x_1, \ldots, x_n]$ such that

$$\varphi(x_i) = \varphi_i(x_1, \ldots, x_n), \text{ where } \varphi_i = \varphi_i(x_1, \ldots, x_n) \in K[x_1, \ldots, x_n],$$

determines a polynomial map $\varphi^* = (\varphi_1, \ldots, \varphi_n): \mathbf{A}^n_K \to \mathbf{A}^n_K$ of the affine space $\mathbf{A}^n_K$ into itself of the form

$$(2.2) \qquad (x_1, \ldots, x_n) \to (\varphi_1(x_1, \ldots, x_n), \ldots, \varphi_n(x_1, \ldots, x_n))$$

The converse is also true: to each polynomial map $\varphi^* : \mathbf{A}_K^n \to \mathbf{A}_K^n$ of the form (2.2) corresponds the above mentioned endomorphism $\varphi$ of the algebra $K[x_1, \ldots, x_n]$. We will make use of this relation below.

Denote by $M_\varphi$ the variety $\varphi^*(\mathbf{A}_K^n)$. We shall say that the variety $M_\varphi$ *corresponds* to the endomorphism $\varphi$ of the polynomial algebra $K[x_1, \ldots, x_n]$. The coordinate ring $K[M_\varphi]$ of the variety $M_\varphi$ is $K[M_\varphi] = K[x_1, \ldots, x_n]/I$, where

$$I = \{H(x_1, \ldots, x_n) \mid H(\varphi_1(\vec{x}), \ldots, \varphi_n(\vec{x})) = 0\}$$

is the ideal in $K[x_1, \ldots, x_n]$ corresponding to the variety $M_\varphi$. It is clear that $K[M_\varphi] \simeq K[\varphi_1(\vec{x}), \ldots, \varphi_n(\vec{x})]$ and $\dim M_\varphi = \operatorname{trdeg} K[\varphi_1(\vec{x}), \ldots, \varphi_n(\vec{x})]$.

**Lemma 2.9.** *The variety $M_\varphi$ is irreducible.*

*Proof.* Since the affine variety $\mathbf{A}_K^n$ corresponding to the algebra $K[x_1, \ldots, x_n]$ is irreducible and the image of an irreducible algebraic variety is also irreducible [6, 21], the variety $M_\varphi$ is irreducible. Hint: coordinate ring of an image isomorhic to subring of the coordinate ring of the preimage, hence has no zero divisors.) □

**Lemma 2.10.** *Let $\phi_1, \phi_2$ be endomorphisms of $K[x_1, \ldots, x_n]$ and $M_{\phi_1}, M_{\phi_2}$ be two corresponding varieties, respectively. The following properties hold:*
  (1) *If $\phi_1 \sim \phi_2$, then $M_{\phi_1} \cong M_{\phi_2}$ and the corresponding coordinate rings are isomorphic.*
  (2) *$\phi_1 \preceq \phi_2$ if and only if the coordinate ring of $M_{\phi_1}$ is a quotient ring of the coordinate ring of $M_{\phi_2}$. In this case $\dim M_{\phi_2} \leq \dim M_{\phi_1}$, where $\dim X$ is the Krull dimension of a variety $X$. If the quotient ring is proper, then the inequality is strict.*

*Proof.*   (1) By item (3) of Lemma 2.5, the coordinate rings of the varieties $M_{\phi_1}$ and $M_{\phi_2}$ are isomorphic. Therefore, the above varieties themselves are isomorphic.
  (2) By item (2) of Lemma 2.5, the coordinate ring of the variety $M_{\phi_1}$ is a quotient ring of the coordinate ring of the variety $M_{\phi_2}$ by some its ideal. As a consequence, $\dim M_{\phi_1} \leq \dim M_{\phi_2}$ (see also [6, 21]).
□

Let $\psi$ be an endomorphism of $K[x_1, \ldots, x_n]$ of "external" rank $m$. The last lemma shows that there exists no chains of endomorphisms $\psi_i$ of the form (2.1) of length more than $m$ beginning with $\psi$. It means that the inner rank of $\psi$ is less or equal than the outer its rank. In order to prove the proposition 2.8 we need to establish an opposite inequality, i.e., to prove that there exists a chain (2.1) of length $m$ beginning with $\psi$.

**Lemma 2.11.** *Notations being as above, let $\dim M_\varphi = m$. Then there exists an endomorphism $\varphi'$ of $K[x_1, \ldots, x_n]$ such that $\varphi' \prec \varphi$ and $\dim M_{\varphi'} = m - 1$.*

The assertion of this lemma is evident for $m = 1$: in this case it is sufficient to consider specialization $x_i \to \xi_i$, $\xi_i \in K$, into ground field $K$.

Now we pass to the general case. We need the following lemma

**Lemma 2.12.** *Let $R$ be a subalgebra of $K[x_1, \ldots, x_n]$ of a transcendence degree $m$ $(m \leq n)$. Then there exists an embedding from $R$ into $K[x_1, \ldots, x_m]$.*

**Remark 2.13.** A similar statement for field embeddings was established in [4].

*Proof.* It is known that any transcendence base of a subalgebra A of a algebra B can be extended to a transcendence base of the algebra B. Let $y_1, \ldots, y_m$ be a transcendence base of $R$. We can complete this base to a base $y_1, \ldots, y_m, z_1, \ldots, z_{n-m}$ of $K[x_1, \ldots, x_n]$. It is clear that the elements $z_1, \ldots, z_{n-m}$ are algebraically independent over $R$ and they generate a subalgebra $R[z_1, \ldots, z_{n-m}]$ of $K[x_1, \ldots, x_n]$. Therefore, the affine domain $R[z_1, \ldots, z_{n-m}]$ can be embedded into an affine domain $K[x_1, \ldots, x_m][x_1, \ldots, x_{n-m}]$. However, it is known that if $A$ and $B$ are two domains such that $A[x_1, \ldots, x_s]$ can be embedded into $B[x_1, \ldots, x_s]$, then $A$ can be embedded into $B$ (see [4]). Therefore, $R$ can be embedded into the polynomial algebra $K[x_1, \ldots, x_m]$. □

Now, by Lemma 2.12 one can assume that polynomials $\varphi_1, \ldots, \varphi_n$ defining the mapping $\varphi$ belong to $K[x_1, \ldots, x_m]$ and $\mathrm{trdeg}(\varphi_1, \ldots, \varphi_n) = m$, $m \leq n$.

**Lemma 2.14.** *Let $\varphi_1(x_1, \ldots, x_m), \ldots, \varphi_n(x_1, \ldots, x_m)$, where $n \geq m$, be a collection of polynomials from $K[x_1, \ldots, x_m]$ which generates the subalgebra of $K[x_1, \ldots, x_n]$ of transcendence degree $m$. Then for any specialization $x_m \to \xi$, $\xi \in K$, except a finite set of values of $\xi \in K$, the algebra $K[\varphi_1(x_1, \ldots, x_{m-1}, \xi), \ldots, \varphi_n(x_1, \ldots, x_{m-1}, \xi)]$ has the transcendence degree $m - 1$.*

*Proof.* Without loss of generality it is sufficient to consider the case when $K$ is an algebraically closed field (tensoring over algebraic closure, if necessary). Consider a mapping $\Phi : \mathbf{A}_K^m \to \mathbf{A}_K^{n+1}$ such that $\Phi(\vec{x}) = (\varphi_1(\vec{x}), \ldots, \varphi_n(\vec{x}), x_m)$ where $\vec{x} = (x_1, \ldots, x_m)$. Denote by $M$ the image of $\Phi$. Since $\mathrm{trdeg}(\varphi_1, \ldots, \varphi_n) = m$ and the dimension of image $\Phi$ is at most $m$, we have $\dim M = m$. Now we consider a projection $\pi : \mathbf{A}_K^{n+1} \to \mathbf{A}_K^1$ such that $\pi(z_1, \ldots, z_n, x_m) = x_m$. Denote by $\pi_1$ the restriction of $\pi$ to $M$. It is clear that $\pi_1$ is an epimorphic mapping. Further we use the following

**Theorem 2.15.** [6, 21] *If $f : X \to Y$ is a regular mapping between irreducible varieties $X$ and $Y$: $f(X) = Y$, $\dim X = n$, $\dim Y = m$, then $m \leq n$ and*

(1) $\dim f^{-1}(y) \geq n - m$ *for every point $y \in Y$.*
(2) *There exists a non empty set $U \subset Y$ such that $\dim f^{-1}(y) = n - m$ for all $y \in U$.*

In our case $Y = \mathbf{A}_K^1$, $\dim Y = 1$, $\dim X = m$. Therefore, for all points of $\mathbf{A}_K^1$, except points of closed subvariety $T$ of $\mathbf{A}_K^1$, the fiber $\pi^{-1}(\xi)$ has the dimension $m - 1$. Therefore,

$$\mathrm{trdeg}\, K[P_1(x_1, \ldots, x_{m-1}, \xi), \ldots, P_n(x_1, \ldots, x_{m-1}, \xi)] = m - 1.$$

except a finite set of $\xi \in K$. This concludes the proof of Lemma 2.14. □

**Remark 2.16.** A proof of Lemma 2.11 follows immediately from the above Lemma in the **case of an infinite ground field**. Indeed, if a field $K$ is infinite, by Lemma 2.14 we can choose $\xi \in K$ such that $\varphi_1' = \varphi_1(x_1, x_2, \ldots, x_{n-1}, \xi), \ldots, \varphi_n' = \varphi_n(x_1, \ldots, x_{n-1}, \xi)$ and $\mathrm{trdeg} K[\varphi_1'(\vec{x}), \ldots, \varphi_n'(\vec{x})] = m - 1$. As a corollary, we have $\dim M_{\varphi'} = k - 1$, where $\varphi' = (\varphi_1', \ldots, \varphi_n')$. Hence, our Lemma 2.11 is proven in the case of an infinite field. This provides a description of the group $\mathrm{Aut}(\mathrm{End}(K[x_1, \ldots, x_n]))$ for the case of an infinite ground field $K$ as was obtained earlier by Berzins [3].

However, in the case of a finite ground field there can be no such small jumps from $\varphi_i$ to $\varphi_i'$, such that $\dim M_{\varphi'} = \dim M_{\varphi} - 1$, for any specialization of variables into a ground field $K$.

**Example 2.17.** Let $|K| = q$ and $\varphi_i = \prod_{k=1}^{n}(x_k^q - x_k) \cdot x_i$. It is evident that $\operatorname{trdeg}(\varphi_1, \ldots, \varphi_n) = n$. However, any specialization of $\varphi_i$ of the form: $x_n \to \xi, \xi \in K$, yields us $\varphi_i' = 0$.

If a field $K$ is finite *instead of specializations of $x_n$ into ground field we consider substitutions into polynomials depending on other variables, in particular, on powers of other variables.* We need the following

**Theorem 2.18.** [4] *Let $\xi_1, \ldots, \xi_s$ be algebraic over $K[x_1, \ldots, x_m]$, the polynomials $Q_i(\vec{t}, \vec{x}, \vec{\xi})$, $i = 1, \ldots, n$, are algebraically independent for some value of set of parameter $\vec{t} = (t_1, \ldots, t_n)$ in some extension field $k_1$ of the ground field $k$. Then there exists polynomials $R_i \in \Phi[x_1]$, $i = 1, 2, \ldots, r$, $\vec{R} = (R_1, \ldots, R_r)$ such that the set of polynomial*

$$\{Q_1(\vec{t}, \vec{x}, \vec{\xi}), \ldots, Q_n(\vec{t}, \vec{x}, \vec{\xi})\}$$

*is algebraically independent. Moreover, if the growth of the sequence*

$$n_1 \ll n_2 \ll \cdots \ll n_r$$

*is sufficiently large, we may be assume $R_i = x_1^{n_i}$. The above statement is still valid if we replace "$k[x_1, \ldots, x_m]$" by "$k(x_1, \ldots, x_m)$" and "polynomial" for rational function. In this case we can put $R_i = x_1^{-n_i}$.*

*Instead of $x_1$ one can take any other variable $x_i$; $\Phi = \mathbb{Z}_p$ if $\operatorname{char} K = p$ and $\Phi = \mathbb{Z}$ if $\operatorname{char} K = 0$.*

We use a special case of this Theorem for $r = 1$ and $s = 0$, i.e, a variant of this Theorem without $\xi_i$. The next Assertion is also needed for the proof of Lemma 2.11 in the case of a finite ground field $K$.

**Assertion 2.19.** Let $Q_1(x_1, \ldots, x_m), \ldots, Q_n(x_1, \ldots, x_m)$ be a set of polynomials from $K[x_1, \ldots, x_m]$, $|K| < \infty$, and the transcendence degree of the algebra

$$K[Q_1(x_1, \ldots, x_m), \ldots, Q_n(x_1, \ldots, x_m)]$$

equal to $m$, where $m > 1$ and $m \leq n$. If $r \in \mathbb{N}$ is sufficiently large, then

$$\operatorname{trdeg}(K[Q_1(x_1, \ldots, x_1^r), \ldots, Q_n(x_1, \ldots, x_1^r)]) = m - 1.$$

*Proof.* Denote by $A = K[Q_1(x_1, \ldots, x_{m-1}, x_1^r), \ldots, Q_n(x_1, \ldots, x_{m-1}, x_1^r)]$. It is clear that $A \subseteq K[x_1, \ldots, x_{m-1}]$, i.e., $\operatorname{trdeg}(A) \leq m - 1$. We have to prove that the opposite inequality is also fulfilled for sufficiently large $r$. Since

$$\operatorname{trdeg}(K[Q_1(x_1, \ldots, x_m), \ldots, Q_n(x_1, \ldots, x_m)]) = m,$$

we can choose $m$ algebraically independent polynomials between $Q_i$. Without loss of generality, we can set that these polynomials are $Q_1, \ldots, Q_m$. By Lemma 2.14, there exists $\eta \in \bar{K}$, where $\bar{K}$ is the algebraic closure of field $K$, such that

$$\operatorname{trdeg}(\bar{K}[Q_1(x_1, \ldots, x_{m-1}, \eta), \ldots, Q_m(x_1, \ldots, x_{m-1}, \eta)]) = m - 1.$$

Without loss of generality, we can suppose that the first $m - 1$ polynomials $Q_i(x_1, \ldots, x_{m-1}, \eta)$, $1 \leq i \leq m - 1$, are algebraically independent over $\bar{K}$. By Theorem 2.18, there exists a natural $r_0$, such that the polynomials

$$Q_1(x_1, \ldots, x_{m-1}, x^r), \ldots, Q_{m-1}(x_1, \ldots, x_{m-1}, x^r)$$

are algebraically independence over $K$ for any $r \geq r_0$. Since the dimension of the subring $K[Q_1(x_1,\ldots,x_{m-1},x^r),\ldots,Q_{m-1}(x_1,\ldots,x_n,x^r)]$ is not less than the dimension of its subring $K[Q_1(x_1,\ldots,x_{m-1},x^r),\ldots,Q_n(x_1,\ldots,x_{m-1},x^r)]$, the proof is complete. $\square$

We summarize our results in the following

**Assertion 2.20.** Let $\varphi = (\varphi_1(x_1,\ldots,x_n),\ldots,\varphi_n(x_1,\ldots,x_n))$ be an endomorphisms of $K[x_1,\ldots,x_n]$ of "internal" rank $m$. Then there exists an endomorphism $\psi = (\psi_1(x_1,\ldots,x_m),\ldots,\psi_n(x_1,\ldots,x_m))$, $\psi_i(x_1,\ldots,x_m) \in K[x_1,\ldots,x_m]$, such that $\varphi \sim \psi$. In addition, an endomorphism

$$\psi'_{(r)} = (\psi_1(x_1,\ldots,x_{m-1},x_1^r),\ldots,\psi_n(x_1,\ldots,x_{m-1},x_1^r))$$

has the rank at most $m-1$ for any $r \in \mathbb{N}$. Moreover, there exists $r_0 \in \mathbb{N}$ such that for all $r \geq r_0$ holds: $\psi'_{(r)} \prec \psi$. As consequence, $\psi'_{(r)} \prec \varphi$ and an "internal" rank of $\psi'_{(r)}$ is equal to $m-1$ for all $r \geq r_0$.

With these Assertion, the proof of Lemma 2.11 is straightforward. Now we ready to prove Proposition 2.8

**Proof of Proposition 2.8** Suppose that $\varphi$ has an "internal" rank $m$, i.e., there exists a maximal chain of length $m$ beginning with $\varphi$:

(2.3) $$\varphi \precsim \varphi_{m-1} \precsim \cdots \precsim \varphi_1 \prec \varphi_0,$$

We have a descending chain of the corresponding varieties $M_{\varphi_i}$:

(2.4) $$M_{\varphi_0} \subseteq M_{\varphi_1} \subseteq \cdots \subseteq M_{\varphi_{m-1}} \subseteq M_\varphi$$

The induction argument on the length $m$ of the chain (2.4) leads us to the case $m = 0$ for which our assertion is evident. Therefore, the "external" rank of $\varphi$ is also equal to $m$.

Conversely, let an endomorphism $\varphi$ be of "external" rank $m$, i.e., trdeg Im $\varphi = m$. By Lemma 2.11, there exists an endomorphism $\psi_{m-1}$ of $K[x_1,\ldots,x_n]$ such that $\psi_{m-1} \prec \varphi$ and $\dim M_{\psi_{m-1}} = m - 1$. In the same way, we can construct a chain of the form (2.3) beginning with $\varphi$. It is clear that this chain has the length $m$, as desired.

Since the chain (2.1) is invariant under automorphisms of End $K[x_1,\ldots,x_n]$, we have

**Corollary 2.21.** Let $\Phi \in \mathrm{Aut}(\mathrm{End}(A))$, $\psi \in \mathrm{End}(A)$, and $\mathrm{rk}\,(\psi) = m$. Then $\mathrm{rk}\,(\Phi(\psi)) = m$.

**Remark 2.22.** Below we need endomorphisms of rank zero and one. By Definition 2.1, an endomorphism $\psi$ of $A$ is of rank zero if $\psi(A) = K$. An endomorphism $\varphi$ of $A$ is of rank one if $\mathrm{trdeg}(\mathrm{Im}\,\varphi) = 1$. It is known [4], [20], that every integrally closed subalgebra $B$ of $A = K[x_1,\ldots,x_n]$ of transcendence degree 1 is isomorphic to a polynomial algebra $K[t]$ in variable $t$. Taking into account that the integer closure $B$ of the algebra $\varphi(A)$ in $A$ is an algebra of the same transcendence degree as $\varphi(A)$, we conclude that the algebra $B$ is isomorphic to a polynomial algebra $K[t]$ in variable $t$. As a consequence, the algebra $\varphi(A)$ is a polynomial algebra $K[y]$, where $y$ is an element in $K[x_1,\ldots,x_n]$.

2.2. **Representations of Kronecker semigroup of rank $n$.** Recall the definition of Kroneker endomorphisms of the free associative algebra $A$.

**Definition 2.23.** (cf. [9, 11]) *Kroneker endomorphisms* of $A$ in the base $X = \{x_1, \ldots, x_n\}$, $x_i \in A$, are the endomorphisms $e_{ij}$, $i, j \in [1n]$, of $A$ which are determined on free generators $x_k \in X$ by the rule: $e_{ij}(x_k) = \delta_{jk} x_i$, $x_i \in X$, $i, j, k \in [1n]$ and $\delta_{jk}$ is the Kronecker delta.

It is clear that any Kronecker endomorphism of $A$ has rank 1.

**Definition 2.24.** A semigroup $\Gamma_n$ with an adjoint zero element 0 generated by $b_{ij}$, $ij \in [1n]$, with defining relations
$$b_{ij} \cdot b_{km} = \delta_{jk} b_{im}, \ b_{ij} \cdot 0 = 0 \cdot b_{ij} = 0$$
is called a *Kronecker semigroup of rank $n$*.

Denote by $E_n$ a semigroup generated by $e_{ij}$, $i, j \in [1n]$, and an adjoint zero. Clearly, the semigroup $E_n$ is a Kronecker semigroup of rank $n$.

**Remark 2.25.** We have a notion of the rank of a Kronecker semigroup $\Gamma$. Don't confuse it with the rank of an endomorphism of $A$.

**Definition 2.26.** *A representation of a semigroup $T$ in the semigroup $\operatorname{End} A$ is a homomorphism $\nu : T \to \operatorname{End} A$.*

**Definition 2.27.** Let $\rho : \Gamma_n \to \operatorname{End} A$ be a representation of the Kronecker semigroup $\Gamma$ of rank $n$ in $\operatorname{End} A$. We say that the representation $\rho$ is *singular* if $\operatorname{rk} \rho(b_{ij}) = 0$ for any $i, j \in [1n]$.

In fact, it is sufficient to require that $\operatorname{rk} \rho(b_{11}) = 0$.

**Proposition 2.28.** *Let $\rho : \Gamma_n \to \operatorname{End} A$ be a singular representation of the Kronecker semigroup $\Gamma$ of rank $n$ in $\operatorname{End} A$ and $q = \rho \cdot \rho^{-1}$ the kernel congruence on $\Gamma_n$. Then $\Gamma_n/q \cong A$, where $A = \langle \varphi \rangle$ is a one-element semigroup such that $\rho(0) = \varphi$, $\varphi \in \operatorname{End} A$, and $\operatorname{rk}(\varphi) = 0$. Conversely, if $\varphi \in \operatorname{End} A$ is an endomorphism of rank 0, then there exists a representation $\rho : \Gamma_n \to \operatorname{End} A$ such that $\rho(0) = \varphi$.*

*Proof.* From $0 \cdot b_{ij} = 0$, $i, j \in [1n]$, it follows $\varphi \rho(b_{ij}) = \varphi$, where $\rho(0) = \varphi$. Since $\varphi$ is the identical mapping on $K$ and $\operatorname{rk}(\rho(b_{ij})) = 0$, we have $\rho(b_{ij}) = \varphi$ for any $i, j \in [1n]$. Thus, $\Gamma_n/q \cong A$, where $A = \langle \varphi \rangle$.

Conversely, if $\varphi$ is an endomorphism of $\operatorname{End} A$ such that $\operatorname{rk}(\varphi) = 0$. Define a representation $\rho : \Gamma_n \to \operatorname{End} A$ by the rule $\rho(0) = \rho(b_{ij}) = \varphi$ for all $i, j \in [1n]$. It is clear that we obtained a required representation $\rho$. $\square$

**Remark 2.29.** Let $\rho : \Gamma_n \to \operatorname{End} A$ be a singular representation of the Kronecker semigroup $\Gamma_n$ of rank $n$ in $\operatorname{End} A$ such that $\rho(0) = \varphi$, $\varphi \in \operatorname{End} A$, and $\operatorname{rk}(\varphi) = 0$. We can set $\varphi(x_i) = \alpha_i$, $\alpha_i \in K$. Denote by $\psi : K^n \to K^n$ the mapping on $K^n$ such that $\psi(x_1, \ldots, x_n) = (x_1 - \alpha_1, \ldots, x_n - \alpha_n)$. Define a representation $\widehat{\rho} : \Gamma_n \to \operatorname{End} A$ of $\Gamma_n$ in $\operatorname{End} A$ by the rule $\widehat{\rho}(0) = \widehat{\rho}(b_{ij}) = \varphi \psi$ for all $i, j \in [1n]$. Then $\varphi \psi = \widehat{O}$ and $\widehat{\rho}(0) = \widehat{O}$.

**Proposition 2.30.** *Let $\rho : \Gamma_n \to \operatorname{End} A$ be a non-singular representation of a Kronecker semigroup $\Gamma_n$. Then, $\operatorname{rk}(\rho(b_{ij})) = 1$ for all $i, j \in [1n]$.*

*Proof.* We will make use below relations between polynomial map $\varphi : K^n \to K^n$ and endomorphisms of the polynomial algebra $K[x_1, \ldots, x_n]$, described on the page 3.

Denote $\rho(b_{ij})$ by $\varphi_{ij}$, $i, j \in [1n]$. Let $\bar{\varphi}_{ij}$ be the endomorphisms of the algebra $B = K[x_1, \ldots, x_n]$ of commutative polynomials in variables $x_1, \ldots, x_n$ induced by the endomorphisms $\varphi_{ij}$ of the algebra $A$. Clearly, $\bar{\varphi}_{ij}\bar{\varphi}_{km} = \delta_{jk}\bar{\varphi}_{im}$. Let us note $o \cdot \varphi_{im} = \widehat{O}$. For a fix $j \in [1n]$ consider $\bar{\varphi}_{jj}$ as a polynomial mapping from $K^n$ into $K^n$, i.e., $\bar{\varphi}_{jj}(x_1, \ldots, x_n) = (\bar{\varphi}_{jj}(x_1), \ldots, \bar{\varphi}_{jj}(x_n))$. Since $\bar{\varphi}_{jj}^2 = \bar{\varphi}_{jj}$, the mapping $\bar{\varphi}_{jj}$ has a fixed point in $K^n$. This point $d = (d_1, \ldots, d_n)$, $d_i \in K$, can be chosen arbitrarily from the image of $\bar{\varphi}_{jj}$. Therefore, we have $\bar{\varphi}_{jj}(d_1, \ldots, d_n) = (d_1, \ldots, d_n)$.

Denote by $T : K^n \to K^n$ the polynomial mapping on $K^n$ such that $T(x_1, \ldots, x_n) = (x_1 + d_1, \ldots, x_n + d_n)$. Let $\tilde{\varphi}_{ij} = T^{-1}\bar{\varphi}_{ij}T$ be a mapping $K^n$ into itself. Denote by $p_{ij}^{(k)}$ the element $T^{-1}\bar{\varphi}_{ij}T(x_k)$. Since the mapping $\tilde{\varphi}_{ii}$ has the fixed point $0 \in K^n$, the elements $p_{ii}^{(k)}$ do not have constant terms for any $i, k \in [1n]$. Now we will prove that the elements $p_{ij}^{(k)}$, $i, j, k \in [1n]$, also do not have constant terms. Assume, on the contrary, that there exist $i, j, k \in [1n]$, $i \neq j$, such that the element $p_{ij}^{(k)}$ has a constant term. Since the elements $p_{jj}^{(m)} = T^{-1}\bar{\varphi}_{jj}T(x_m)$ do not have a constant term for any $m, j \in [1n]$, we obtain

$$(T^{-1}\bar{\varphi}_{jj}T)(T^{-1}\bar{\varphi}_{ij}T)(x_k) = (T^{-1}\bar{\varphi}_{jj}T)p_{ij}^{(k)} \neq 0.$$

On the other hand, since $i \neq j$

$$(T^{-1}\bar{\varphi}_{jj}T)(T^{-1}\bar{\varphi}_{ij}T)(x_k) = (T^{-1}\bar{\varphi}_{jj}\bar{\varphi}_{ij}T)(x_k) = 0.$$

This contradiction proves that the elements $p_{ij}^{(k)} = T^{-1}\bar{\varphi}_{ij}T(x_k)$ do not have a constant term for any $i, j, k \in [1n]$. As a consequence, the elements $T^{-1}\varphi_{ij}T(x_k)$ do not have constant terms for any $i, j, k \in [1n]$, too.

Denote the mapping $T^{-1}\varphi_{ij}T : A \to A$ by $\hat{\varphi}_{ij}$. We now prove that $\hat{\varphi}_{ij}(A)$ is a subalgebra of $K[w]$ for some $w \in A$. Let $I$ be the ideal of $A$ generated by $x_1, \ldots, x_n$. Since the elements $\hat{\varphi}_{ij}(x_k)$, $i, j, k \in [1n]$, do not have a constant term, $\hat{\varphi}_{ij}(I^s) \subseteq I^s$ for any $s \geq 1$. Now we fix some $i, j \in [1n]$ and consider induced maps $\tilde{\varphi}_{ij}^{(s)} : I^s/I^{s+1} \to I^s/I^{s+1}$ for any $s \geq 1$. We intend to prove that $\operatorname{Im}\tilde{\varphi}_{ij}^{(s)}$ are one-dimensional vector spaces over $K$. Let $s = 1$. Then $\tilde{\varphi}_{ij}^{(1)} : I/I^2 \to I/I^2$ is a linear mapping from the vector space $I/I^2$ into itself. Since $\tilde{\varphi}_{ij}^{(1)}\tilde{\varphi}_{mk}^{(1)} = \delta_{jm}\tilde{\varphi}_{ik}^{(1)}$, by Lemma 4.7 [11] there exists a basis $\bar{z}_{r1} = z_r + I^2$, where $z_r \in I$, $r \in [1n]$, of $I/I^2$ such that $\tilde{\varphi}_{ij}^{(1)}(\bar{z}_{r1}) = \delta_{jr}\bar{z}_{i1}$. For a fix number $s \geq 2$ denote $\bar{z}_{rs} = z_r + I^{s+1}, r \in [1n]$. We have $\tilde{\varphi}_{ij}^{(s)}(\bar{z}_{i_1 s} \cdots \bar{z}_{i_s s}) = \delta_{ji_1} \cdots \delta_{ji_s}\bar{z}_{is}^s$. Thus, $\tilde{\varphi}_{ij}^{(s)}(I^s/I^{s+1})$ is a one-dimensional vector space with a basis $\{\bar{z}_{is}^s\}$. The latter assertion holds for any $s \geq 2$. As a consequence, we have $\hat{\varphi}_{ij}(A) \subseteq K[z_i]$. Hence, $\varphi_{ij}(A)$ is a subalgebra of $K[w]$, where $w = Tz_i$. Since the representation $\rho$ of $\Gamma$ is non-singular, $K \subset \varphi_{ij}(A)$. Thus, $\operatorname{rk}(\varphi_{ij}) = \operatorname{rk}\rho(b_{ij}) = 1$ for all $i, j \in [1n]$. $\square$

### 2.3. Bases and subbases of the semigroup $\operatorname{End} A$. We need the following

**Definition 2.31.** A set of endomorphisms $\mathcal{B}_e = \{e'_{ij} | e'_{ij} \in \operatorname{End} A \text{ and } e'_{ij} \neq \widehat{O}, \forall i, j \in [1n]\}$ of $A$ is called a *subbase* of $\operatorname{End} A$ if $e'_{ij}e'_{km} = \delta_{jk}e'_{im}, \forall i, j, k, m \in [1n]$.

Let us note that $0 \cdot e'_{ij} = \widehat{O}$. Denote by $E'$ a semigroup of $\operatorname{End} A$ generated by endomorphisms $e'_{ij}$ and the endomorphism $\widehat{O}$. By Theorem 2.30, we obtain the following

**Corollary 2.32.** $\operatorname{rk}(e'_{ij}) = 1$ for any $i, j \in [1n]$.

We can assume that $e'_{ij}(A)$ is a subalgebra of $K[z_{ij}]$, $i, j \in [1n]$, where $z_{ij} \in A$. For the sake of simplicity we write $z_{ii} = z_i$, $i \in [1n]$.

**Definition 2.33.** ("External" definition of a base collection of $\operatorname{End} A$.) We say that the subbase $\mathcal{B}_e$ is a *base collection of endomorphisms* of $A$ (or *a base* of $\operatorname{End} A$, for short) if $Z = \{z_i \mid z_i \in A \text{ such that } e'_{ii}(A) \subseteq K[z_i], i \in [1n]\}$ is a base of $A$.

Now we show that there exists a subbase of $\operatorname{End} A$ that is not its base.

**Example 2.34.** Let $\varphi_{ij} : K[x_1, x_2] \to K[x_1, x_2]$, where $i, j \in \{1, 2\}$, be endomorphisms of the free associative-commutative algebra $A = K[x_1, x_2]$ such that

(2.5)  $\varphi_{11}(x_1) = x_1 + x_1 x_2,\ \varphi_{11}(x_2) = 0,\ \varphi_{22}(x_1) = 0,\ \varphi_{22}(x_2) = x_2,$
$\varphi_{12}(x_1) = 0,\ \varphi_{12}(x_2) = x_1 + x_1 x_2,\ \varphi_{21}(x_1) = x_1,\ \varphi_{21}(x_2) = 0.$

It is easy to see that $\operatorname{rk}(\varphi_{ij}) = 1$ and $\varphi_{ij} \varphi_{km} = \delta_{jk} \varphi_{im}$ for any $i, j, k, m \in \{1, 2\}$, i.e., the set of endomorphisms $B_\varphi = \{\varphi_{ij} | \varphi_{ij} \in \operatorname{End} A,\ i, j \in \{1, 2\}\}$ is a subbase of the semigroup $\operatorname{End} A$. We will prove that $B_\varphi$ is not its base. It is clear that $\varphi_{11}(A) = K[u]$, where $u = x_1 + x_1 x_2$, and $\varphi_{22}(A) = K[x_1]$. We can take $z_1 = u$ and $z_2 = x_1$. The elements $z_1$ and $z_2$ generate the algebra $K[x_1 + x_1 x_2, x_1]$. Let us show that $K[x_1 + x_1 x_2, x_2] \neq K[x_1, x_2]$. If, on the contrary, $K[x_1 + x_1, x_2, x_2] = K[x_1, x_2]$ then $x_1 = \alpha(x_1 + x_1 x_2) + \beta x_2 + P(u, x_2)$, where $\deg P(u, x_2) \geq 2$ and $\alpha, \beta \in K$. Hence $\beta = 0, \alpha = 1$ and $P(u, x_2) = 0$. We come to a contradiction. Therefore, the subbase $B_\varphi$ is not a base of $\operatorname{End} A$.

"Internal" definition of a *base collection* of $\operatorname{End} A$ is a bit tricky (see [11, 9]). It was inspired by G.Zhitomirski (see [22]).

**Definition 2.35.** ("Internal" definition of a base collection of $\operatorname{End} A$.) The subbase of endomorphisms $\mathcal{B}_e = \{e'_{ij} | e'_{ij} \in \operatorname{End} A,\ i, j \in [1n]\}$ of $\operatorname{End} A$ is its base if for any collection of endomorphisms $\alpha_i : A \to A,\ \forall i \in [1n]$, and any subbase $\mathcal{B}_f = \{f'_{ij} \mid i, j \in [1n]\}$ of $\operatorname{End} A$ there exist endomorphisms $\varphi, \psi \in \operatorname{End} A$ such that

(2.6) $$\alpha_i \circ f'_{ii} = \psi \circ e'_{ii} \circ \varphi,\ \text{for all } i \in [1n].$$

Our aim is to prove statement similar to the proposition 2.27 in [5].

**Proposition 2.36.** *Internal and external definitions of a base collection of $\operatorname{End} A$ are equivalent.*

*Proof.* Let a subbase of endomorphisms $\mathcal{B}_e$ be a base according Definition 2.33. Since $\operatorname{rk}(f'_{ij}) = 1,\ \forall i, j \in [1n]$, there exist elements $y_{ij} \in A,\ i, j \in [1n]$, such that $K \subset f'_{ij}(A(X)) \subseteq K[y_{ij}]$ for all $i, j \in [1n]$. Define endomorphisms $\psi$ and $\varphi$ of $A$ as follows:

$$\varphi(x_i) = z_i \text{ and } \psi(z_i) = \alpha_i(y_i),\ \text{for all } i \in [1n],$$

where $e'_{ii}(A) \subseteq K[z_i],\ z_i \in A$, and $y_i = y_{ii},\ \forall i \in [1n]$. Since $Z = \langle z_i \mid z_i \in A, i \in [1n]\rangle$ is a base of $A$, the endomorphism $\psi$ is well-defined. Now it is easy to check that the condition (2.6) with the given $\varphi$ and $\psi$ is fulfilled.

Conversely, assume that the condition (2.6) is fulfilled for the subbase $\mathcal{B}_e$. Let us prove that $Z = \langle z_i \mid z_i \in A, i \in [1n]\rangle$ is a base of $A$. Choosing $\alpha_i = e_{ii}$ and $f'_{ij} = e_{ij}$, $i,j \in [1n]$, in (2.6), we obtain

$$e_{ii} = \psi \circ e'_{ii} \circ \varphi,$$

i.e., $\psi(e'_{ii}\varphi(x_i)) = x_i$ for any $i \in [1n]$. Denote by $t_i = e'_{ii}\varphi(x_i)$. We have $\psi(t_i) = x_i$. Since $A$ is Hopfian, i.e., any surjective endomorphism of A into itself is isomorphism, the elements $t_i$, $i \in [1n]$, form the base of $A$. By Corollary and Remark 2.22 2.32, $K \subset e'_{ii}(A) \subseteq K[z_i]$. Therefore, there exists a non-scalar polynomial $\chi_i(z_i) \in K[z_i]$ such that $t_i = \chi_i(z_i)$. Since $t_i = \chi_i(z_i), i = 1,\ldots, n$, forms the base of $A$, the elements $z_i, i = 1,\ldots, n$, forms a base of $A$ as claimed. $\square$

Now we deduce

**Corollary 2.37.** *Let $\Phi \in \operatorname{Aut}\operatorname{End} A$ and $E$ be the subsemigroup of $\operatorname{End} A$ generated by the Kronecker endomorphisms $e_{ij}$, $i,j \in [1n]$ (see Definition 2.23). Then $\mathcal{C} = \{\Phi(e_{ij}) \mid i,j \in [1n]\}$ is a base of $\operatorname{End} A$.*

*Proof.* Assume that $\operatorname{rk}(\Phi(e_{ij})) = 0$ for some $i,j \in [1n]$. By Corollary 2.21, we obtain $\operatorname{rk}(e_{ij}) = 0$. We arrived at a contradiction. Thus, $\operatorname{rk}(\Phi(e_{ij})) \neq 0$. Since $\Phi(e_{ij})\Phi(e_{km}) = \delta_{jk}\Phi(e_{im})$, the set $\mathcal{C}$ is a subbase of $\operatorname{End} A$. It is easy to check that the condition (2.6) is fulfilled for the subbase $\mathcal{C}$. Thus, $\mathcal{C}$ is a base of $\operatorname{End} A$. $\square$

**Lemma 2.38.** *Let $\mathcal{B}_e = \{e'_{ij} \mid e'_{ij} \in \operatorname{End} A, i,j \in [1n]\}$ be a base collection of endomorphisms of $\operatorname{End} A$. Then there exists a base $Z' = \{z'_k \mid z'_k \in A, k \in [1n]\}$ of $A$ such that the endomorphisms $e'_{ij}$ from $\mathcal{B}_e$ are Kronecker ones of $A$ in $Z'$.*

*Proof.* With the preceding notation from Definition 2.33 we have that the equality $(e'_{ii})^2 = e'_{ii}$ implies $e'_{ii}(z_i) = z_i$, $i \in [1n]$. Since $e'_{ii}e'_{ij}(z_j) = e'_{ij}(z_j)$ and $K \subset e'_{ii}(A) \subseteq K[z_i]$, there exists a non-scalar polynomial $f_j(z_i) \in K[z_i]$ such that $e'_{ij}(z_j) = f_j(z_i)$. Similarly, there exists a non-scalar polynomial $g_i(z_j) \in K[z_j]$ such that $e'_{ji}(z_i) = g_i(z_j)$. We have

$$z_j = e'_{jj}(z_j) = e'_{ji}e'_{ij}(z_j) = e'_{ji}(f_j(z_i)) = f_j(g_i(z_j)) \text{ for all } i,j \in [1n]$$

and, in a similar way, $z_i = g_i(f_j(z_i))$ for all $i,j \in [1n]$. Thus $f_j$ and $g_i$ are linear polynomials over $K$ in variables $z_i$ and $z_j$, respectively. Therefore,

$$(2.7) \qquad e'_{ij}(z_j) = a_i z_i + b_i, \ a_i, b_i \in K \text{ and } a_i \neq 0.$$

Note that $e'_{ij}(z_k) = e'_{ij}(e'_{kk}(z_k)) = 0$ if $k \neq j$. Now we have for $i \neq j$

$$0 = {e'_{ij}}^2(z_j) = e'_{ij}(a_i z_i + b_i) = e'_{ij}(b_i) = b_i,$$

i.e., $e'_{ij}(z_j) = a_i z_i$, $a_i \neq 0$. Let $z'_i = a_i^{-1} z_i$. We obtain a base $Z = \{z'_k \mid z'_k \in A, k \in [1n]\}$ of $A$ such that $e'_{ij}(z'_k) = \delta_{jk} z'_k$, $i,j,k \in [1n]$, i.e., $e'_{ij}$ are Kronecker endomorphisms of $A$ in the base $Z'$. The proof is completed. $\square$

## 3. Automorphisms of the semigroup $\operatorname{End} A$

3.1. **On the group** $\operatorname{Aut}\operatorname{End} A$. We need the following notion.

**Definition 3.1.** ([7]) Let $A_1$ and $A_2$ be algebras over $K$ from a variety $\mathcal{A}$, $\delta$ be an automorphism of $K$ and $\varphi : A_1 \to A_2$ be a ring homomorphism of these algebras. A pair $(\delta, \varphi)$ is called a *semi-linear homomorphism* from $A_1$ to $A_2$ if
$$\varphi(\alpha \cdot u) = \delta(\alpha) \cdot \varphi(u), \ \forall \alpha \in K, \ \forall u \in A_1.$$

**Definition 3.2.** [17] An automorphism $\Phi$ of the semigroup $\operatorname{End} A$ of endomorphisms of $A$ is called *quasi-inner* if there exists an *adjoined bijection* $s : A \to A$ such that $\Phi(\nu) = s\nu s^{-1}$, for any $\nu \in \operatorname{End} A$

**Definition 3.3.** [17] A quasi-inner automorphism $\Phi$ of $\operatorname{End} A$ is called *semi-inner* if there exists a field automorphism $\delta : K \to K$ such that $(\delta, s)$ is a semi-linear automorphism of $A$, i.e., for any $\alpha \in K$ and $a, b \in A$ the following conditions hold:
1. $s(a + b) = s(a) + s(b)$,
2. $s(a \cdot b) = s(a) \cdot s(b)$,
3. $s(\alpha a) = \delta(\alpha) s(a)$.

We say that the pair $(\delta, s)$ defines the semi-inner automorphism $\Phi$ of $A$ with the *adjoined ring automorphism s*. If $\delta$ is the identity automorphism of $K$, we call the automorphism $\Phi$ *inner*.

The description of quasi-inner automorphisms of $\operatorname{End} A$ is as follows.

**Proposition 3.4.** [3, 9, 11] *Let $\Phi \in \operatorname{Aut} \operatorname{End} A$ be a quasi-inner automorphism of $\operatorname{End} A$. Then $\Phi$ is of semi-inner automorphisms of $\operatorname{End} A$.*

We will use the following fact:

**Proposition 3.5.** [9, 11] *Let $\Phi \in \operatorname{Aut} \operatorname{End} A$ and $E$ be the subsemigroup of $\operatorname{End} A$ generated by $e_{ij}$, $i, j \in [1n]$. Elements of the semigroup $\Phi(E)$ are Kronecker endomorphisms of $A$ in some base $U = \{u_1, \ldots, u_n\}$, $u_i \in A$, if and only if $\Phi$ is a quasi-inner automorphism of $\operatorname{End} A$.*

Now we obtain one of the main result of the paper

**Theorem 3.6.** *Every automorphism of the group $\operatorname{Aut} \operatorname{End} A$ is semi-inner.*

*Proof.* By Corollary 2.37, the set of endomorphisms $\mathcal{C} = \{\Phi(e_{ij}) \mid \forall i \in [1n]\}$ is a base collection of endomorphisms of $A$. By Lemma 2.38, there exists a base $S = \langle s_k \mid s_k \in A, k \in [1n] \rangle$ such that the endomorphisms $\Phi(e_{ij})$ are Kronecker endomorphisms in $S$. According to Proposition 3.5, we obtain that $\Phi$ is quasi-inner. By virtue of Proposition 3.4, every automorphism the group $\operatorname{Aut} \operatorname{End} A$ is semi-inner and as claimed. □

**Remark 3.7.** If $\mathcal{CA}$ is the category of commutative-associative algebras over a field $K$, we take $\mathcal{SCA}$ to be the category with objects all associative algebras from the category $\mathcal{A}$, morphisms all pairs $\psi_\delta = (\psi, \delta) : A \to B$, $A, B \in \operatorname{Ob} \mathcal{SA}$, such that $\psi : A \to B$ are ring homomorphisms from $A$ to $B$, $\delta : K \to K$ are automorphisms of the field $K$ and $\psi_\delta(\lambda a) = \lambda^\delta \psi(a)$, $a \in A$. Morphisms $\psi_\delta$ of the category $\mathcal{SA}$ are called *semi-linear homomorphisms*(or *semihomomorphisms*) from $A$ to $B$ (cf. Definition 3.1). Denote by $\operatorname{SEnd} A$ the semigroup of semiendomorphisms of $A$ with the usual composition of maps in the category $\mathcal{SCA}$.

Clearly, that the definitions of endomorphisms of rank one and zero can be transfer to the category $\mathcal{SCA}$. All results about bases and subbases from the sections 2.3 are also true. As a consequence, we obtain the following

**Theorem 3.8.** *Every automorphism of the group $\operatorname{Aut} \operatorname{SEnd} A$ is semi-inner.*

## 4. Automorphisms of the category $\mathcal{A}^\circ$

Recall the following notions of the category isomorphism and equivalence (cf. [12]). *An isomorphism $\varphi : \mathcal{C} \to \mathcal{M}$ of categories* is a functor $\varphi$ from $\mathcal{C}$ to $\mathcal{M}$, which is a bijection both on objects and morphisms. In other words, there exists a functor $\psi : \mathcal{M} \to \mathcal{C}$ such that $\psi\varphi = 1_\mathcal{C}$ and $\varphi\psi = 1_\mathcal{M}$.

Let $\varphi_1$ and $\varphi_2$ be two functors from $\mathcal{C}_1$ to $\mathcal{C}_2$. *A functor isomorphism* $s : \varphi_1 \longrightarrow \varphi_2$ is a collection of isomorphisms $s_D : \varphi_1(D) \longrightarrow \varphi_2(D)$ defined for all $D \in \mathrm{Ob}\,\mathcal{C}_1$ such that for every $\nu : D \longrightarrow B$, $\nu \in \mathrm{Mor}\,\mathcal{C}_1$, $B \in \mathrm{Ob}\,\mathcal{C}_1$

$$s_B \cdot \varphi_1(\nu) = \varphi_2(\nu) \cdot s_D$$

holds, i.e., the following diagram

$$\begin{array}{ccc} \varphi_1(D) & \xrightarrow{s_D} & \varphi_2(D) \\ \varphi_1(\nu) \downarrow & & \downarrow \varphi_2(\nu) \\ \varphi_1(B) & \xrightarrow{s_B} & \varphi_2(B) \end{array}$$

is commutative. An isomorphism of functors $\varphi_1$ and $\varphi_2$ is denoted by $\varphi_1 \cong \varphi_2$.

*An equivalence of categories* $\mathcal{C}$ and $\mathcal{M}$ is a pair of functors $\varphi : \mathcal{C} \to \mathcal{M}$ and $\psi : \mathcal{M} \to \mathcal{C}$ such that $\psi\varphi \cong 1_\mathcal{C}$ and $\varphi\psi \cong 1_\mathcal{M}$. If $\mathcal{C} = \mathcal{M}$, then we get the notions of *automorphism* and *autoequivalence* of the category $\mathcal{C}$.

For every small category $\mathcal{C}$, denote the group of all its automorphisms by $\mathrm{Aut}\,\mathcal{C}$. We distinguish the following classes of automorphisms of $\mathcal{C}$.

**Definition 4.1.** [8, 15] An automorphism $\varphi : \mathcal{C} \to \mathcal{C}$ is *equinumerous* if $\varphi(D) \cong D$ for any object $D \in \mathrm{Ob}\,\mathcal{C}$ ; $\varphi$ is *stable* if $\varphi(D) = D$ for any object $D \in \mathrm{Ob}\,\mathcal{C}$ ; and $\varphi$ is *inner* if $\varphi$ and $1_\mathcal{C}$ are naturally isomorphic, i.e., $\varphi \cong 1_\mathcal{C}$.

In other words, an automorphism $\varphi$ is inner if for all $D \in \mathrm{Ob}\,\mathcal{C}$ there exists an isomorphism $s_D : A \to \varphi(D)$ such that

$$\varphi(\nu) = s_B \nu s_D^{-1} : \varphi(D) \to \varphi(B)$$

for any morphism $\nu \in \mathrm{Mor}_\mathcal{C}(A, B)$.

Denote by $\mathrm{EqnAut}\,\mathcal{C}$, $\mathrm{StAut}\,\mathcal{C}$, and $\mathrm{Int}\,\mathcal{C}$ the collections of equinumerous, stable, and inner automorphisms of the group $\mathrm{Aut}\,\mathcal{C}$, respectively.

Let $\Theta$ be a variety of linear algebras over $K$. Denote by $\Theta^0$ *the full subcategory of finitely generated free algebras* $F(X), |X| < \infty$, of the variety $\Theta$. Consider a constant morphism $\nu_0 : F(X) \to F(X)$ such that $\nu_0(x) = x_0$, $x_0 \in F(X)$, for every $x \in X$.

**Theorem 4.2.** *(Reduction Theorem [8, 13, 16, 22]) Let the free algebra $F(X)$ generate a variety $\Theta$, and $\varphi \in \mathrm{StAut}\,\Theta^0$. If $\varphi$ acts trivially on the monoid $\mathrm{Mor}_{\Theta^0}(F(X), F(X))$ and $\varphi(\nu_0) = \nu_0$, then $\varphi$ is inner, i.e., $\varphi \in \mathrm{Int}\,\Theta^0$.*

Define the notion of a semi-inner automorphism of the category $\Theta^0$ of free finitely generated algebras in the category $\Theta$.

**Definition 4.3.** [15] An automorphism $\varphi \in \mathrm{Aut}\,\Theta^0$ is called *semi-inner* if there exists a family of semi-isomorphisms $\{s_{F(X)} = (\delta, \tilde{\varphi}) : F(X) \to \tilde{\varphi}(F(X)), F(X) \in \mathrm{Ob}\,\Theta^0\}$, where $\delta \in \mathrm{Aut}\,K$ and $\tilde{\varphi}$ is a ring isomorphism from $F(X)$ to $\tilde{\varphi}(F(X))$ such that for any homomorphism $\nu : F(X) \longrightarrow F(Y)$ the following diagram

$$F(X) \xrightarrow{s_{F(X)}} \tilde{\varphi}(F(X))$$
$$\nu \downarrow \qquad \downarrow \varphi(\nu)$$
$$F(Y) \xrightarrow{s_{F(Y)}} \tilde{\varphi}(F(Y))$$

is commutative.

Further, we will need the following

**Proposition 4.4.** [8, 15] *For any equinumerous automorphism $\varphi \in \operatorname{Aut}\mathcal{C}$ there exist a stable automorphism $\varphi_S$ and an inner automorphism $\varphi_I$ of the category $\mathcal{C}$ such that $\varphi = \varphi_S \varphi_I$.*

Now we give a description of the groups $\operatorname{Aut}\mathcal{CA}^\circ$ over any field. Note that a description of this group over infinite fields was given in [2]

**Theorem 4.5.** *All automorphisms of the group $\operatorname{Aut}\mathcal{A}^\circ$ of automorphisms of the category $\mathcal{CA}^\circ$ are semi-inner automorphisms of the category $\mathcal{CA}^\circ$.*

*Proof.* Let $\varphi \in \operatorname{Aut}\mathcal{A}^\circ$. It is clear that $\varphi$ is an equinumerous automorphism. By Proposition 4.4, $\varphi$ can be represented as a composition of a stable automorphism $\varphi_S$ and an inner automorphism $\varphi_I$. Since stable automorphisms does not change free algebras from $\mathcal{A}^\circ$, we obtain that $\varphi_S \in \operatorname{Aut}\operatorname{End}A$. By Theorem 3.6, $\varphi_S$ is semi-inner of End $A$. Using this fact and Reduction Theorem 4.2, we obtain that all automorphisms of the group $\operatorname{Aut}\mathcal{CA}^\circ$ are semi-inner automorphisms of the category $\mathcal{CA}^\circ$. This completes the proof. □

**Problem 4.6.** *Describe the groups $\operatorname{Aut}\mathcal{B}^\circ$ and $\operatorname{Aut}\operatorname{End}B$, where $B = B(x_1, \ldots, x_n)$, is a free algebra of a non-associative variety $\mathcal{B}$ of linear algebras finitely generated by a set $X = \{x_1, \ldots, x_n\}$.*

Note that the above mentioned groups were described for some homogeneous varieties of linear algebras in [5, 9, 11]. In particular, a description of these group for the variety of all Lie algebras over any field was obtained there. A corresponding description in the case of Lie algebras over any infinite field was obtained in [15, 22].

## 5. Acknowledgments

The authors are grateful to B.Plotkin for attracting their attention to this problem and interest to this work. The first author was supported by the Israel Science Foundation (grant No. 1178/06).

[1]*Department of Mathematics*, Department of Mathematics, Bar Ilan University, Ramat Gan, 52900, Israel
   *E-mail address*: `belova@macs.biu.ac.il`

[2]*Department of Mathematics*, Ben Gurion University, Beer Sheva, 84105, Israel
   *E-mail address*: `lipyansk@cs.bgu.ac.il`